\documentclass[11pt]{article}
\usepackage{CJK}
\usepackage{amsmath}
\usepackage{amsmath,amssymb,latexsym,color}
\usepackage[mathscr]{eucal}
\usepackage{graphicx}
\usepackage{pstricks}
\newtheorem{thm}{Theorem}[section]
\newtheorem{defin}[thm]{Definition}
\newtheorem{prop}[thm]{Proposition}
\newtheorem{lemma}[thm]{Lemma}
\newtheorem{cor}[thm]{Corollary}

\newtheorem{example}{Example}

\newcommand{\remark}{\vspace{2ex}\noindent{\bf Remark.\quad}}

\newcommand{\proof}{{\it Proof.\quad}}
\newcommand{\qed}{\hfill\Box\medskip}
\begin{document}
\begin{CJK*}{GBK}{song}

\title{\bf Identifying codes of lexicographic product of graphs}

\author{Min Feng\quad Min Xu\quad Kaishun Wang\footnote{Corresponding author. E-mail address: wangks@bnu.edu.cn}\\
{\footnotesize   \em  Sch. Math. Sci. {\rm \&} Lab. Math. Com. Sys.,
Beijing Normal University, Beijing, 100875,  China}}
\date{}
\maketitle

\begin{abstract}
Gravier et  al. \cite{Gra} investigated  the identifying codes of
Cartesian product of two graphs. In this  paper we consider the
identifying codes of lexicographic product $G[H]$ of a connected
graph $G$ and an arbitrary graph $H$, and obtain the minimum
cardinality of identifying codes of $G[H]$ in terms of some
parameters of $G$ and $H$.

\medskip
\noindent {\em Key words:} Identifying code; lexicographic product.

\end{abstract}
\bigskip

\bigskip

\section{Introduction}
In this paper, we only consider finite undirected simple graphs with
at least two vertices. For a given graph $G$, we often write $V(G)$
for the vertex set of $G$ and $E(G)$ for the edge set of $G$. For
any two vertices $u$ and $v$ of $G$, $d_{G}(u,v)$ denotes the
distance between $u$ and $v$ in $G$. Given a vertex $v\in V(G)$, we
define $B_{G}(v)=\{u|u\in V(G),d_G(u,v)\leq 1\}$. A \emph{code} $C$
is a nonempty set of vertices. For a code $C$, we say that $C$
\emph{covers} $v$ if $B_{G}(v)\cap C\neq\emptyset$; We say that $C$
\emph{separates} two distinct vertices $x$ and $y$ if $B_{G}(x)\cap
C\neq B_{G}(y)\cap C$. An \emph{identifying code} of $G$ is a code
which covers all the vertices of $G$ and separates any pair of
distinct vertices of $G$.
 If $G$ admits at least one identifying code, we say $G$ is \emph{identifiable} and denote
 the minimum cardinality of all identifying codes of $G$ by $I(G)$.

The concept of identifying codes was  introduced by Karpovsky et al.
\cite{Ka1} to model a fault-detection problem in multiprocessor
systems.  It was noted in \cite{Ch,Coh} that determining the
identifying code with the minimum cardinality in a graph is an
NP-complete problem. Many researchers  have focused on the study of
identifying codes in some restricted classes of graphs, for example,
paths \cite{Be}, cycles \cite{Be,Gr,Xu}, and hypercubes
\cite{Bl1,Ho,Ka2,Mo}.

Gravier et  al. \cite{Gra} investigated  the identifying codes of
Cartesian product of two cliques. In this  paper, we consider the
identifying codes of lexicographic product $G[H]$ of a connected
graph $G$ and an arbitrary graph $H$. In Section 2, we introduce two
new families of codes which are closely related to identifying
codes, and compute the minimum cardinalities of the two codes for
paths and cycles, respectively. In Section 3, we give the sufficient
and necessary condition when $G[H]$ is identifiable, and obtain the
minimum cardinality of identifying codes of $G[H]$ in terms of some
parameters of $G$ and $H$.

\section{Two new families of codes}

For  a graph $H$, let $C'\subseteq V(H)$ be a code which separates
any pair of distinct vertices of $H$, we use $I'(H)$ to denote the
minimum cardinality of all possible $C'$; let $C''\subseteq V(H)$ be
a code which separates any pair of distinct vertices of $H$ and
satisfies $C''\not\subseteq B_{H}(v)$ for every $v\in V(H)$, we use
$I''(H)$ to denote the minimum cardinality of all possible $C''$.

The two parameters $I'(H)$ and $I''(H)$ are used to compute the
minimum cardinality of identifying codes of $G[H]$ of graphs $G$ and
$H$  (see Theorem~\ref{main}). In this section we shall compute the
two parameters for  paths and  cycles, respectively.

Given an integer $n\geq 3$, let $P_{n}$ be the path of order $n$ and
$C_{n}$ be the cycle of order $n$. Suppose
 $$\begin{array}{c}
V(P_{n})=\{0,1,\ldots,n-1\},\;
 E(P_{n})=\{ij|j=i+1,i=0,\ldots,n-2\};\\
 V(C_{n})=\mathbb{Z}_{n}=\{0,1,\ldots,n-1\},\; E(C_{n})=\{ij|j=i+1,i\in \mathbb{Z}_{n}\}.
\end{array}$$

\begin{example}\rm  $I'(P_{3})=2$ and $I''(P_{3})$ is not well defined; $I'(P_{4})=3$
and $I''(P_{4})=4$; $I'(P_{5})=I''(P_{5})=3$; $I'(P_{6})=3$ and
$I''(P_{6})=4$.

For $P_{4}$, $\{0,1,2\}$ is an identifying code, but
$\{0,1,2\}\subseteq B_{P_{4}}(1)$ and $\{0,1,3\}$ can not separate
$0$ and $1$. For $P_{5}$, $\{0,2,4\}$ separates any pair of distinct
vertices. For $P_{6}$, $\{1,2,3\}$ separates any pair of distinct
vertices, but  $\{1,2,3\}\subseteq B_{P_{6}}(2)$.
\end{example}

\begin{example}\rm   $I'(C_{4})=3$ and $I''(C_{4})=4$; $I'(C_{5})=3$ and
$I''(C_{5})=4$; $I'(C_{6})=I''(C_{6})=3$; $I'(C_{7})=I''(C_{7})=4$;
$I'(C_{9})=I''(C_{9})=6$; $I'(C_{11})=I''(C_{11})=6$.

For $C_{4}$, $\{0,1,2\}$ is an identifying code, but
$\{0,1,2\}\subseteq B_{C_{4}}(1)$. For $C_{5}$, $\{0,1,2\}$ is
an identifying code, but $\{0,1,2\}\subseteq B_{C_{5}}(1)$ and $\{0,1,3\}$ can not separate $0$ and $1$.
For $C_{6}$, both $\{3,4,5\}$ and $\{0,2,4\}$ separate any pair of
distinct vertices. For $C_{7}$, $\{3,4,5,6\}$ separates any pair of
distinct vertices. For $C_{9}$,  both $\{3,4,5,6,7,8\}$ and
$\{0,2,4,6,7,8\}$ separate  any pair of distinct vertices. For
$C_{11}$, $\{3,4,5,8,9,10\}$ separates any pair of distinct
vertices.
\end{example}

The minimum cardinality of identifying codes of a path or a cycle
was computed in \cite{Be, Gr}.

\begin{prop}\label{pc} {\rm (\cite{Be, Gr})} {\rm(i)} For $n\geq3$,
$I(P_{n})=\lfloor\frac{n}{2}\rfloor+1$;

{\rm(ii)} For $n\geq 6$, $I(C_{n})=\left\{
\begin{array}{ll}
\frac{n}{2}, &n ~\textup{is even},\\
\frac{n+3}{2}, &n ~\textup{is odd}.
\end{array}\right.$
\end{prop}

In order to  compute the two parameters  for  paths  and cycle,
 we need the following useful lemma.

\begin{lemma}\label{rk'}
Let $H$ be an identifiable graph.

{\rm (i)} $I(H)-1\leq I'(H)\leq I(H)$;

{\rm (ii)} If $\Delta(H)\leq|V(H)|-2$, then $I(H)-1\leq I'(H)\leq I''(H)\leq
I(H)+1$, where $\Delta(H)$ is the maximum degree of   $H$.
\end{lemma}
\proof Let $C'$ be a code which separates any pair of distinct
vertices of $H$.

(i) Since there exists at most one vertex $v$ not covered by $C'$,
  $C'\cup \{v\}$  is an identifying code of $H$.

(ii) Note that there exists at most one vertex $v$ such that
$C'\subseteq B_H(v)$. Since $\Delta(H)\leq|V(H)|-2$, there exists
$v_{0}\in V(H)\backslash B_H(v)$ such that $C''=C'\cup\{v_{0}\}$  is
a code which separates any pair of distinct vertices of $H$ and
satisfies $C''\not\subseteq B_{H}(w)$ for every $w\in V(H)$. It
follows that $I'(H)\leq I''(H)\leq I'(H)+1$. By (i), (ii) holds.
$\qed$

For two integers   $i\leq j$, let $[i,j]=\{i,i+1,\ldots,j\}$.

\begin{prop}\label{p}
For $n\geq 7$, $I'(P_{n})=I''(P_{n})=\lfloor\frac{n}{2}\rfloor+1$.
\end{prop}
\proof By Lemma \ref{rk'}, $I'(P_{n})=I(P_{n})$ or $I(P_{n})-1$. If
$I'(P_{n})=I(P_{n})-1$, then there exists a code $W'$ of size
$I(P_{n})-1$ such that $W'$ separates any pair of distinct vertices
of $P_{n}$ and $B_{P_{n}}(i_{0})\cap W'=\emptyset$ for a unique
$i_{0}\in[0,n-1]$.

\emph{Case 1}. $n$ is odd. Let $W=([0,i_{0}]\cap
W')\cup\{i-1|i\in[i_{0}+1,n-1]\cap W'\}\subseteq[0,n-2]$.  Since $W$ covers all
vertices of $P_{n-1}$, $W$ is an identifying code of $P_{n-1}$. By
Proposition \ref{pc},
$$\frac{n+1}{2}=I(P_{n-1})\leq|W|=|W'|=I(P_{n})-1=\frac{n-1}{2},$$ a contradiction.

\emph{Case 2}. $n$ is even.  By Proposition \ref{pc},
$|W'|=I(P_{n})-1=\frac{n}{2}$.

\emph{Case 2.1}. $i_{0}\neq 0$ and $i_{0}\neq n-1$. Then
$i_{0}-1,i_{0},i_{0}+1\not\in W'$,
and $i_{0}-2,i_{0}-3,i_{0}-4,i_{0}+2,i_{0}+3,i_{0}+4\in W'$,
so $4\leq i_{0}\leq n-5$. Let $W=W'\cap[0,i_{0}-1]$ and
$\overline{W}=\{i-i_{0}-1|i\in W'\cap [i_{0}+1,n-1]\}$. Then $W$ is an identifying code of $P_{i_{0}}$ and $\overline{W}$ is
an identifying code of $P_{n-i_{0}-1}$.  By Proposition \ref{pc}, we have
$$\frac{n}{2}=|W'|=|W|+|\overline{W}|\geq
I(P_{i_{0}})+I(P_{n-i_{0}-1})=\lfloor\frac{i_{0}}{2}\rfloor+1+\lfloor\frac{n-i_{0}-1}{2}\rfloor+1=\frac{n+2}{2},$$
 a contradiction.

\emph{Case 2.2}. $i_{0}=0$ or $n-1$. Without loss of generality, assume
$i_{0}=n-1$. Then $n-1,n-2\not\in W'$, and
$n-3,n-4,n-5\in W'$. We can observe the following results:
\begin{align}
\label{equa1}&|W'\cap[i,i+3]|\geq 2, i\in[0,n-4],\\
\label{equa2}&|W'\cap[0,2]|\geq2,\\
\label{equa3}&|W'\cap[0,4]|\geq3.
\end{align}

\emph{Case 2.2.1}. $n=4k$. By (\ref{equa1}) and (\ref{equa2}),
$2k=|W'|\geq 2\lfloor\frac{n-5-3}{4}\rfloor+3+2=2k+1$, a
contradiction.

\emph{Case 2.2.2}. $n=4k+2$. By (\ref{equa1}) and (\ref{equa3}),
$2k+1=|W'|\geq 2\lfloor\frac{n-5-5}{4}\rfloor+3+3=2k+2$, a
contradiction.

Therefore, $I'(P_{n})=I(P_n)$. Note that $I''(P_{n})=I'(P_{n})$ when
$I'(P_{n})\geq 4$. By Proposition \ref{pc}, the desired result
follows. $\qed$

\begin{prop}\label{c}
$I'(C_{n})=I''(C_{n})=\left\{
\begin{array}{ll}
\frac{n}{2}, &n \textup{ is even and }n\geq 8,\\
\frac{n+3}{2}, &n \textup{ is odd and }n\geq 13.
\end{array}\right.$
\end{prop}
\proof If $I'(C_{n})<\lceil\frac{n}{2}\rceil$, then there exists a
code $W'$ such that $|W'|=I'(C_{n})<\lceil\frac{n}{2}\rceil$ and
$W'$ separates any pair of distinct vertices of $C_{n}$. It follows
that there exists $i_{0}$ such that $i_{0}, i_{0}+1\not\in W'$.
Without loss of generality, assume $n-1, 0\not\in W'$. Since $W'$ is
also a subset of $V(P_{n})$ and $B_{C_{n}}(j)\cap
W'=B_{P_{n}}(j)\cap W'$ for any $j\in[0,n-1]$, $W'$ separates any
pair of distinct vertices of $P_{n}$. By Proposition \ref{p},
$\lceil\frac{n}{2}\rceil\leq
I'(P_{n})\leq|W'|<\lceil\frac{n}{2}\rceil$, a contradiction. Hence
$I'(C_n)\geq\lceil\frac{n}{2}\rceil$.

\emph{Case 1}. $n$ is even and $n\geq 8$. By Proposition \ref{pc}
and Lemma \ref{rk'}, $\frac{n}{2}\leq I'(C_{n})\leq
I(C_{n})=\frac{n}{2}$. Hence $I'(C_{n})=\frac{n}{2}$.

\emph{Case 2}. $n$ is odd and $n\geq 13$.  By Proposition \ref{pc}
and Lemma \ref{rk'}, $I'(C_{n})=\frac{n+3}{2}$ or $\frac{n+1}{2}$.
If $I'(C_{n})=\frac{n+1}{2}$, then there exists a code $W'$ of
size $\frac{n+1}{2}$ such that $W'$ separates any pair of
distinct vertices of $C_{n}$ and $B_{C_{n}}(i_{0})\cap W'=\emptyset$
for a unique $i_{0}\in[0,n-1]$. Without loss of generality, assume
$i_{0}=1$.  Then $0,1,2\not\in W'$ and $3,4,5, n-3, n-2, n-1\in W'$.
We can observe the following results:
\begin{align}
\label{equa4}&|W'\cap[i,i+3]|\geq 2,i\in[6,n-7],\\
\label{equa5}&|W'\cap[6,11]|\geq 3.
\end{align}

\emph{Case 2.2.1}. $n=4k+1$. By (\ref{equa4}), $2k+1=|W'|\geq
2\lfloor\frac{n-9}{4}\rfloor+6=2k+2$, a contradiction.

\emph{Case 2.2.2}. $n=4k+3$. By (\ref{equa4}) and (\ref{equa5}),
$2k+2=|W'|\geq 2\lfloor\frac{n-9-6}{4}\rfloor+6+3=2k+3$, a
contradiction.

Therefore, $I'(C_{n})=\frac{n+3}{2}$.

Since $I''(C_{n})=I'(C_{n})$ when $I'(C_{n})\geq 4$, the desired
result follows.$\qed$

\section{Main results}

The \emph{lexicographic product} $G[H]$ of graphs $G$ and $H$  is
the graph with the vertex set $V(G)\times V(H)=\{(u,v)|u\in V(G),v\in
V(H)\}$, and the edge set $\{\{(u_1,v_1),(u_2,v_2)\}|\\ d_G(u_1,u_2) =1, \textup{ or } u_1=u_2\textup{ and }d_H(v_1,v_2)=1\}$.
 For any two distinct vertices $(u_{1},v_{1})$, $(u_{2},v_{2})$ of $G[H]$, we observe that
\begin{equation}\label{distance}
d_{G[H]}((u_{1},v_{1}),(u_{2},v_{2}))=\left\{
\begin{array}{ll}
1,                 &\textup{if}~u_{1}=u_{2},d_{H}(v_{1},v_{2})=1,\\
2,                 &\textup{if}~u_{1}=u_{2},d_{H}(v_{1},v_{2})\geq2,\\
d_{G}(u_{1},u_{2}),&\textup{if}~u_{1}\neq u_{2}.
\end{array}\right.
\end{equation}

For $u\in V(G)$, let $N_{G}(u)=B_{G}(u)\backslash\{u\}$. For any
$u_{1},u_{2}\in V(G)$, define $u_{1}\equiv u_{2}$ if and only if
$B_{G}(u_{1})=B_{G}(u_{2})$ or $N_{G}(u_{1})=N_{G}(u_{2})$.
Hernando et al. \cite{He}  proved that $``\equiv"$ is an equivalent
relation and the equivalence class of a vertex is of three types: a
class of size $1$, a clique of size at least $2$, an independent set
of size at least $2$. Denote all equivalence classes by
\begin{equation}\label{class}
W_{1},\ldots,W_{p},U_{1},\ldots,U_{k},V_{1},\ldots,V_{l},
\end{equation}
where

(i) $|W_{q}|=1$, $q=1,\ldots,p$;

(ii) for any $u_{1},u_{2}\in U_{i},i=1,\ldots,k$,
$B_{G}(u_{1})=B_{G}(u_{2})$;

(iii) for any $u_{1},u_{2}\in V_{j},j=1,\ldots,l$,
$N_{G}(u_{1})=N_{G}(u_{2})$. \\* Denote
$s(G)=|U_{1}|+\cdots+|U_{k}|-k$, $t(G)=|V_{1}|+\cdots+|V_{l}|-l$. We
give an algorithm of computing $s(G)$ and $t(G)$ in Appendix.

For $u\in V(G)$ and $C\subseteq V(H)$, let $C^{u}=\{(u,v)|(u,v)\in V(G[H]), v\in C\}$.  For $S\subseteq V(G[H])$, let $S_{u}=\{v|v\in V(H), (u,v)\in S\}$. Note that $(S_u)^u=H^u\cap S$, where $H^{u}=(V(H))^{u}$.
 By (\ref{distance}), we have
\begin{equation}\label{equ}
B_{G[H]}((u,v))=(B_{H}(v))^{u}\cup\bigcup_{w\in N_{G}(u)} H^{w},
\end{equation}
\begin{equation}\label{equa}
B_{G[H]}((u,v))\cap S=((B_{H}(v))\cap S_u)^u\cup\bigcup_{w\in N_{G}(u)}(S_w)^w.
\end{equation}

In the rest of this section we always assume that $G$ is a connected
graph and $H$ is an arbitrary graph.

\begin{thm}\label{exi}
The lexicographic product $G[H]$ of graphs $G$ and $H$ is identifiable if
and only if

{\rm (i)}  $H$ is identifiable and $\Delta(H)\leq |V(H)|-2$, or

{\rm (ii)} both $G$ and $H$ are identifiable.
\end{thm}

\proof Suppose $G[H]$ is identifiable.  If $H$ is not identifiable,
then there exist two distinct vertices $v_{1},v_{2}$ of $H$ with
$B_{H}(v_{1})=B_{H}(v_{2})$. By (\ref{equ}),
$B_{G[H]}((u,v_{1}))=B_{G[H]}((u,v_{2}))$ for $u\in V(G)$. This
contradicts the condition that $G[H]$ is identifiable.

If $\Delta(H)=|V(H)|-1$ and $G$ is not identifiable, then
there exist  $v\in V(H)$  and two distinct vertices $u_{1},u_{2}$ of $G$ such that
$$B_{H}(v)=V(H) \textup{ and }B_{G}(u_{1})=B_{G}(u_{2}).$$
By (\ref{equ}),  we have
$$
B_{G[H]}((u_{1},v))=H^{u_{1}}\cup\bigcup_{u\in N_{G}(u_{1})} H^{u}=\bigcup_{u\in B_{G}(u_{1})} H^{u}=\bigcup_{u\in B_{G}(u_{2})} H^{u}=B_{G[H]}((u_{2},v)).
$$
This contradicts the condition that $G[H]$ is
identifiable.

Therefore, (i) or (ii) holds.

Conversely, suppose (i) or (ii) holds.  Assume that $G[H]$ is not
identifiable. Therefore, there exist two distinct  vertices
$(u_{1},v_{1}),(u_{2},v_{2})$  such that
$B_{G[H]}((u_{1},v_{1}))=B_{G[H]}((u_{2},v_{2}))$. If $u_{1}\neq
u_{2}$, then $d_G(u_1,u_2)=1$. It follows that
$B_{G}(u_{1})=B_{G}(u_{2})$ and $B_{H}(v_{1})=B_{H}(v_{2})=V(H)$,
contrary to (i) and (ii). If $u_{1}=u_{2}$, then $v_{1}\neq v_{2}$.
By (\ref{equ}), $B_{H}(v_{1})=B_{H}(v_{2})$, contrary to the
condition that $H$ is identifiable. $\qed$

\remark  Let $r$ be a positive integer and $\Gamma$ be a graph. Given a vertex $v\in V(\Gamma)$, define $B_{\Gamma}^{(r)}(v)=\{u|u\in V(\Gamma),d_\Gamma(u,v)\leq r\}$. An $r$-\emph{identifying code} of $\Gamma$ is a code which $r$-covers all the vertices of $\Gamma$ and $r$-separates any pair of distinct vertices of $\Gamma$ (see \cite{Ka1} for details). Identifying codes in this paper are $1$-identifying codes. If $r\geq 2$, then $G[H]$ does not admit any $r$-identifying code. Indeed, by (\ref{distance}), $B_{G[H]}^{(r)}((u,v_1))=B_{G[H]}^{(r)}((u,v_2))$ for $r\geq 2$.

\begin{lemma}\label{geq}
 If $S$ is an identifying code of $G[H]$, then for any vertex $u$ of
$G$, $S_{u}$ separates any pair of distinct vertices
of $H$. Moreover, with reference to (\ref{class}),

{\rm (i)} if $k\neq 0$, then there exists at most one vertex $u\in
U_{i}$ satisfying   $S_{u}\subseteq B_{H}(v)$ for a vertex $v$ of
$H$, where $i=1,\ldots,k$;

{\rm (ii)} if $l\neq 0$, then there exists at most one vertex $u\in
V_{j}$ satisfying   $S_{u}\cap B_{H}(v)=\emptyset$ for a vertex $v$
of $H$, where $j=1,\ldots,l$.
\end{lemma}
\proof  Assume that there exist $u_{0}\in V(G)$ and two distinct vertices $v_{1},v_{2}$ of $H$ such that $S_{u_{0}}\cap B_{H}(v_{1})=S_{u_{0}}\cap B_{H}(v_{2})$.
By (\ref{equa}), $B_{G[H]}((u_{0},v_{1}))\cap S=B_{G[H]}((u_{0},v_{2}))\cap S$, 
contrary to the condition that $S$ is an identifying code of $G[H]$.

(i) Assume that there exist two distinct vertices $u_{1},u_{2}\in
U_{i}$ such that $S_{u_{1}}\subseteq B_{H}(v_{1})$ and
$S_{u_{2}}\subseteq B_{H}(v_{2})$. Since
$B_{G}(u_{1})=B_{G}(u_{2})$, by (\ref{equa})  we have
\begin{equation*}
B_{G[H]}((u_{1},v_{1}))\cap S=(S_{u_{1}})^{u_{1}}\cup\bigcup_{u\in N_{G}(u_{1})}(S_u)^u
                             =\bigcup_{u\in B_{G}(u_{2})}(S_u)^u
                             =B_{G[H]}((u_{2},v_{2}))\cap S.
\end{equation*}
Since $S$ is an identifying code of $G[H]$, $(u_{1},v_{1})=(u_{2},v_{2})$, a contradiction.

(ii) Assume that there exist two different vertices $u_{1},u_{2}\in
V_{j}$ such that $S_{u_{1}}\cap B_{H}(v_{1})=S_{u_{2}}\cap
B_{H}(v_{2})=\emptyset$. Since $N_{G}(u_{1})=N_{G}(u_{2})$, by
(\ref{equa})  we have
\begin{equation*}
B_{G[H]}((u_{1},v_{1}))\cap S=\bigcup_{u\in N_{G}(u_{1})}(S_u)^u=\bigcup_{u\in N_{G}(u_{2})}(S_u)^u
                            =B_{G[H]}((u_{2},v_{2}))\cap S.
\end{equation*}
Since $S$ is an identifying code of $G[H]$, $(u_{1},v_{1})=(u_{2},v_{2})$, a contradiction. $\qed$

In equivalence classes (\ref{class}) of $V(G)$,
choose $\overline{u}_{i}\in U_{i},i=1,\ldots,k$,
and $\overline{v}_{j}\in V_{j},j=1,\ldots,l$.
Let
$\overline{W}_{0}=\cup_{q=1}^{p}W_{q}\cup\{\overline{u}_{1},\ldots,\overline{u}_{k},\overline{v}_{1},\ldots,\overline{v}_{l}\}$
and
$\overline{U}_{i}=U_{i}\backslash\{\overline{u}_{i}\},i=1,\ldots,k$, $\overline{V}_{j}=V_{j}\backslash\{\overline{v}_{j}\},j=1,\ldots,l$. Therefore, we have a partition of $V(G)$:
\begin{equation}\label{partition}
\overline{W}_{0},\overline{U}_{1},\ldots,\overline{U}_{k},\overline{V}_{1},\ldots,\overline{V}_{l}.
\end{equation}

\begin{lemma}\label{leq}
Let $C$ be an identifying code of graph $H$, and let $C',C''$ be two
codes which separate any pair of distinct vertices of $H$ and
$C''\not\subseteq B_{H}(v)$ for every vertex $v$ of $H$. With
reference to (\ref{partition}),
\begin{equation*}
S=\bigcup_{u\in\overline{W}_{0}}(C')^{u}\cup\bigcup_{i=1}^{k}\bigcup_{u\in\overline{U}_{i}}(C'')^{u}\cup\bigcup_{i=1}^{l}\bigcup_{u\in\overline{V}_{i}}C^{u}
\end{equation*}
is an identifying code of $G[H]$.
\end{lemma}
\proof For any $u\in V(G)$, we have
\begin{equation*}
S_{u}=\left\{
\begin{array}{ll}
C', &\textup{ if }u\in \overline{W}_{0},\\
C'', &\textup{ if }u\in \cup_{i=1}^{k}\overline{U}_{i},\\
C,&\textup{ if }u\in \cup_{j=1}^{l}\overline{V}_{j}.
\end{array}\right.
\end{equation*}
Since $G$ is connected, there exists a vertex $w$ adjacent to $u$. By (\ref{distance}), $S$ covers all vertices of $G[H]$.
For any two distinct vertices
$(u_{1},v_{1}),(u_{2},v_{2})\in V(G[H])$, we only need to show that
\begin{equation}\label{inequ}
B_{G[H]}((u_{1},v_{1}))\cap S\neq B_{G[H]}((u_{2},v_{2}))\cap S.
\end{equation}
To prove (\ref{inequ}), it is sufficient to show that there exists
$(u_{0},v_{0})\in S$ such that
\begin{equation}\label{in1}
d_{G[H]}((u_{0},v_{0}),(u_{1},v_{1}))\leq 1,~d_{G[H]}((u_{0},v_{0}),(u_{2},v_{2}))\geq 2
\end{equation}
or
\begin{equation}\label{in2}
d_{G[H]}((u_{0},v_{0}),(u_{2},v_{2}))\leq 1,~d_{G[H]}((u_{0},v_{0}),(u_{1},v_{1}))\geq 2.
\end{equation}

\emph{Case 1}. $u_{1}\not\equiv u_{2}$. Then there exists $u_{0}\in
V(G)\backslash\{u_{1},u_{2}\}$ such that $d_{G}(u_{1},u_{0})=1$ and
$d_{G}(u_{2},u_{0})\geq 2$,  or $d_{G}(u_{1},u_{0})\geq 2$ and
$d_{G}(u_{2},u_{0})=1$. Take $v_0\in S_{u_0}$. Then
$(u_{0},v_{0})\in S$. By (\ref{distance}), (\ref{in1}) or
(\ref{in2}) holds.

\emph{Case 2}. $u_{1}\equiv u_{2}$.

\emph{Case 2.1}. $u_{1}=u_{2}$. Since $S_{u_1}$ separates $v_1$ and $v_2$, $B_{H}(v_{1})\cap S_{u_1}\neq B_{H}(v_{2})\cap S_{u_1}=B_{H}(v_{2})\cap S_{u_2}$. By (\ref{equa}), (\ref{inequ}) holds.

\emph{Case 2.2}. $u_{1}\neq u_{2}$ and $B_{G}(u_{1})=B_{G}(u_{2})$.
Then $u_{1}$ and $u_{2}$ are adjacent and fall into some $U_i$. It
follows that $u_1\in \overline{U}_i$ or $u_2\in \overline{U}_i$.
Without loss of generality, suppose  $u_{1}\in \overline{U}_i$. Pick
$u_{0}=u_{1}$. Since $C''\not\subseteq B_{H}(v_{1})$, there exists
$v_{0}\in C''$ such that $(u_{0},v_{0})\in S$ and
$d_{H}(v_{0},v_{1})\geq 2$. By (\ref{distance}), (\ref{in2}) holds.

\emph{Case 2.3}. $u_{1}\neq u_{2}$ and $N_{G}(u_{1})=N_{G}(u_{2})$.
Then $u_{1}$ and $u_{2}$ are at distance $2$ and fall into some
$V_j$. It follows that $u_1\in \overline{V}_j$ or $u_2\in
\overline{V}_j$. Without loss of generality, suppose  $u_{1}\in
\overline{V}_j$. Pick $u_{0}=u_{1}$. Since $C$ covers $v_{1}$, there
exists $v_{0}\in C$ such that $(u_{0},v_{0})\in S$ and
$d_{H}(v_{0},v_{1})\leq 1$. By (\ref{distance}), (\ref{in1}) holds.
$\qed$

\begin{thm}\label{main}  Suppose {\rm (i)} or {\rm (ii)} holds in Theorem \ref{exi}.

{\rm (i)} If $\Delta(H)\leq|V(H)|-2$, then
\begin{equation}\label{equ2}
I(G[H])=(|V(G)|-s(G)-t(G))I'(H)+s(G)I''(H)+t(G)I(H);
\end{equation}

{\rm (ii)} If $\Delta(H)=|V(H)|-1$, then
\begin{equation}\label{equ1}
I(G[H])=(|V(G)|-t(G))I'(H)+t(G)I(H).
\end{equation}
\end{thm}

\proof
 (i) By Theorem \ref{exi}, $I(H)$ and $I'(H)$ are well defined. Since
$V(H)$ separates any pair of distinct vertices of $H$ and
$V(H)\not\subseteq B_{H}(v)$ for every $v\in V(H)$, $I''(H)$ is well
defined.

Let $S$ be an identifying code of $G[H]$ with the minimum cardinality, by Lemma \ref{geq},
\begin{equation*}
\begin{array}{lll}
I(G[H])&=&|S|=\sum_{i=1}^{p}\sum_{u\in W_{i}}|S_{u}|+\sum_{i=1}^{k}\sum_{u\in U_{i}}|S_{u}|+\sum_{i=1}^{l}\sum_{u\in V_{i}}|S_{u}|\\
&\geq&(p+k+l)I'(H)+(\sum_{i=1}^{k}|U_{i}|-l)I''(H)+(\sum_{i=1}^{l}|V_{i}|-l)I(H)\\
&=&(|V(G)|-s(G)-t(G))I'(H)+s(G)I''(H)+t(G)I(H).
\end{array}
\end{equation*}

Let $C$ be an identifying code of $H$ with the minimum cardinality.
Let $C'$ and $C''$ be two codes with the minimum cardinality such
that they separate any pair of distinct vertices of $H$ and
$C''\not\subseteq B_{H}(v)$ for every vertex $v$ of $H$. By Lemma
\ref{leq},
\begin{equation*}
I(G[H])\leq|S|=(|V(G)|-s(G)-t(G))I'(H)+s(G)I''(H)+t(G)I(H).
\end{equation*}

Therefore, (\ref{equ2}) holds.

(ii) By Theorem \ref{exi}, both $G$ and $H$ are identifiable.
So $I(H)$ and $I'(H)$ are well defined.
Owing to $B_{G}(u_{1})\neq B_{G}(u_{2})$ for any two distinct vertices $u_{1},u_{2}$ of $G$,
we get $k=0$ in (\ref{class}) and (\ref{partition}).
Similar to the proof of (i), (\ref{equ1}) holds. $\qed$

Combining   Propositions \ref{pc}, \ref{p}, \ref{c} and Theorem
\ref{main}, we have

\begin{cor}
Let $G$ be a connected graph of order $m$ $(m\geq 2)$.

$(\textup{i})$ For $n\geq 7$,
$I(G[P_{n}])=m(\lfloor\frac{n}{2}\rfloor+1)$;

$(\textup{ii})$ For $n\geq 12$, $I(G[C_{n}])=\left\{
\begin{array}{ll}
\frac{mn}{2}, &n ~\textup{is even},\\
\frac{m(n+3)}{2}, &n ~\textup{is odd}.
\end{array}\right.$
\end{cor}

\section*{Appendix}

\scriptsize\begin{tabular}[t]{rl} \hline
\quad&\textbf{Algorithm} \\
\hline
\textbf{Input}  & Graph $G$ \\
\textbf{Output} & $W_{1},\ldots,W_{p},U_{1},\ldots,U_{k},V_{1},\ldots,V_{l}$  //the equivalent classes of $V(G)$\\
                & s(G),t(G)\\
\hline \textbf{Step 1}. &Preparation//Input the adjacent matrix $A$ of $G$ and $A+E$ ($E$ is an identity matrix).\\
1.&$V(G)=\{1,\ldots,m\}$;$E(G)=\{ij|ij$ are adjacent in $G$ $\}$\\
2.&\textbf{for} $i=1,\ldots,m$ \textbf{do}\\
3.&\quad \textbf{for} $j=1,\ldots,m$ \textbf{do}\\
4.&\quad\quad \textbf{if} $j=i$ \textbf{then} $a_{ij}:=0$ and $\overline{a}_{ij}:=1$\\
5.&\quad\quad \textbf{else if} $ij\in E$ \textbf{then} $a_{ij}:=1$ and $\overline{a}_{ij}:=1$\\
6.&\quad\quad\quad\quad~\textbf{else} $a_{ij}:=0$ and
$\overline{a}_{ij}:=0$\\
7.&\quad\quad\quad \quad~\textbf{end-if} \\
8.&\quad\quad \textbf{end-if}\\
9.&\quad \textbf{end-for} \\
10.&\textbf{end-for} \\
11.&\textbf{for} $i=1,\ldots,m$ \textbf{do}\\
12.&$A_{i}:=(a_{i1},\ldots,a_{im})$;
$\overline{A}_{i}:=(\overline{a}_{i1},\ldots,\overline{a}_{im})$ \\
13.&\textbf{end-for}\\
\textbf{Step 2}.&Output the equivalent classes of $V(G)$\\
14.&$i:=1$;~$p:=1$;~$k:=1$;~$l:=1$;~$I:=\emptyset$\\
15.&\textbf{while} $i\leq m$ \textbf{do}\\
16.&\quad \textbf{if} $i\in I$ \textbf{then} $i:=i+1$ //$i\equiv i_{0}$ for some $i_{0}<i$\\
17.&\quad \textbf{else if} $i\leq m-1$ \textbf{then}  $W_{p}=\{i\}$;~$U_{k}=\{i\}$;~$V_{t}=\{i\}$ and \textbf{do}\\
18.&\quad\quad\quad\quad\textbf{for} $j=i+1,\ldots,m$ \textbf{do}\\
19.&\quad\quad\quad \quad\quad\textbf{if} $\overline{A}_{j}=\overline{A}_{i}  $\textbf{then} $I:=I\cup\{j\}$ and $U_{k}:=U_{k}\cup\{j\}$ //$B_{G}(j)=B_{G}(i)$\\
20.&\quad\quad\quad \quad\quad\textbf{else if} $A_{j}=A_{i}$ \textbf{then} $I:=I\cup\{j\}$ and $V_{l}:=V_{l}\cup\{j\}$ //$N_{G}(j)=N_{G}(i)$\\
21.&\quad\quad\quad \quad\quad\textbf{end-if}\\
22.&\quad\quad \quad\quad\textbf{end-for}\\
23.&\quad\quad \quad\quad\textbf{if} $|U_{k}|>1$ \textbf{then} \textbf{output} $U_{k}$ and $k:=k+1$\\
24.&\quad\quad \quad\quad\textbf{else if} $|V_{l}|>1$ \textbf{then} \textbf{output} $V_{l}$ and $l:=l+1$\\
25.&\quad\quad \quad\quad\quad\quad~\textbf{else} \textbf{output} $W_{p}$ and $p:=p+1$  //$i\not\equiv j$ for any $j\in V(G)$\\
26.&\quad\quad \quad\quad\quad\quad~\textbf{end-if}\\
27.&\quad\quad \quad\quad\textbf{end-if}\\
28.&\quad\quad\quad\quad $i:=i+1$\\
29.&\quad\quad\quad~\textbf{else} $W_{p}:=\{i\}$, $i:=i+1$ and \textbf{output} $W_{p}$ \\
30.&\quad \quad\quad~\textbf{end-if}\\
31.&\quad \textbf{end-if}\\
32.&\textbf{end-while}\\
\textbf{Step 3}.&Compute $s(G)$ and $t(G)$\\
33.&$s=0$; $t=0$\\
34.&\textbf{If} $k>1$ \textbf{then}\\
35.&\quad \textbf{for} $i=1,\ldots,k-1$ \textbf{do}\\
36.&\quad\quad $s:=s+|U_{i}|$\\
37.&\quad \textbf{end-for}\\
38.&\textbf{If} $k>1$ \textbf{then}\\
39.&\quad \textbf{for} $i=1,\ldots,k-1$ \textbf{do}\\
40.&\quad\quad $t:=t+|V_{i}|$\\
41.&\quad \textbf{end-for}\\
42.& \textbf{output} $s(G)=s$ and $t(G)=t$\\
\hline
\end{tabular}\normalsize

\section*{Acknowledgement} This research is supported by NSF of China (10871027), NCET-08-0052,  and   the
Fundamental Research Funds for the Central Universities of China.

\end{CJK*}
\end{document}